\documentclass[12pt,draft]{amsart}

\theoremstyle{definition}

\RequirePackage{amsmath}
\RequirePackage{amssymb}
\RequirePackage{amsthm}
 \RequirePackage{epsfig}

\begin{document}

\date{}

\title[On the Borel-Cantelli lemma]
{On the M\'{o}ri-Sz\'{e}kely conjectures for the Borel-Cantelli
lemma}

\author{Chunrong Feng and Liangpan Li}

\address{Department of Mathematical Sciences,
Loughborough University,  LE11 3TU,  UK}
 \email{C.Feng@lboro.ac.uk, L.Li@lboro.ac.uk}



\date{}

\begin{abstract}
The purpose of this note is to show by constructing counterexamples
that two conjectures of M\'{o}ri and Sz\'{e}kely for the
Borel-Cantelli lemma are false.
\end{abstract}

\maketitle

\section{Introduction}
Let $\{A_n\}_{n=1}^{\infty}$ be an arbitrary sequence of events in a
probability space $(\Omega,\mathbb{P})$ and denote by $A_{\infty}$
the event that infinitely many $A_n$ occurs simultaneously, i.e.
\[A_{\infty}=\bigcap_{n=1}^{\infty}\bigcup_{k=n}^{\infty}A_k.\]
The classical Borel-Cantelli lemma states that:  if
$\sum_{n=1}^{\infty}\mathbb{P}(A_n)<\infty$ then
$\mathbb{P}(A_{\infty})=0$; else if
$\sum_{n=1}^{\infty}\mathbb{P}(A_n)=\infty$ and
$\{A_n\}_{n=1}^{\infty}$ are mutually independent, then
$\mathbb{P}(A_{\infty})=1$. In the past century many investigations
were devoted to the second implication in the attempt to weaken the
independence condition on $\{A_n\}_{n=1}^{\infty}$. For example, one
of the most applicable results is due to Erd\H{o}s and R\'{e}nyi
(\cite{ErdosRenyi,Renyi}, see also
\cite{Dawson,FengLiShen,KochenStone,Spitzer}) who proved that if
$\sum_{n=1}^{\infty}\mathbb{P}(A_n)=\infty$, then
\[
\mathbb{P}(A_{\infty})\geq\limsup_{n\rightarrow\infty}
\frac{\displaystyle\big(\sum_{k=1}^n{\mathbb{P}}(A_k)\big)^2}{\displaystyle\sum_{i=1}^n\sum_{j=1}^n{\mathbb{P}}(A_iA_j)}=
\limsup_{n\rightarrow\infty}\frac{\displaystyle1}{\displaystyle\mathbb{E}(\alpha_n^2)}\doteq\mbox{ER}(\{A_n\}),
\]
where $\mathbb{E}$ denotes the expectation function,
$\alpha_n\doteq{\displaystyle(\sum_{i=1}^n\mathbb{I}_{A_i})}/{\displaystyle(\sum_{i=1}^n\mathbb{P}({A_i}))}$,
and $\mathbb{I}_{A_i}$ is the indicator function of $A_i$. Later on,
by studying convex and concave Young functions M\'{o}ri and
Sz\'{e}kely (\cite{Mori}, see also
\cite{Amghibech,Frolov,Hu,Liu,Xie})
 improved the Erd\H{o}s-R\'{e}nyi bound to
\[
\mathbb{P}(A_{\infty})\geq\sup_{p\in(0,\infty)\backslash\{1\}}\big(
\limsup_{n\rightarrow\infty}(\mathbb{E}(\alpha_n^p)^{\frac{1}{1-p}})\big)=\lim_{p\searrow0}
\big(\limsup_{n\rightarrow\infty}\mathbb{E}(\alpha_n^p)\big)\doteq
\mbox{MS}(\{A_n\}).
\]
They also proposed the following two conjectures:

\textsc{Conjecture 1}:
\[\mathbb{P}(A_{\infty})=\sup\limits_{\tau:\mathbb{N}\rightarrow\mathbb{N}\ \mbox{is increasing}
}\mbox{MS}(\{A_{\tau(n)}\}).\]

\textsc{Conjecture 2}: If we have an estimate of the form
$\mathbb{P}(A_{\infty})\geq L_k$ $(k\geq2)$ where the constant $L_k$
depends only on $P(A_{i_1}A_{i_2}\cdots A_{i_k})$, $1\leq i_1 \leq
i_2\leq \cdots\leq i_k$, then
\[\sup\limits_{\tau:\mathbb{N}\rightarrow\mathbb{N}\ \mbox{is increasing}}\mbox{ER}(\{A_{\tau(n)}\})\geq L_k.\]

The purpose of this note is to show that both conjectures are false.

\section{A counterexample to Conjecture 1}

Let  $\mathbb{P}$ be the Lebesgue measure on $\Omega=[0,1]$ and let
\[A_{2^{i-1}+k}=[\frac{k}{2^i},\frac{k+1}{2^i}]\cup[\frac{1}{2},1]\ \
\ \ \ i\in\mathbb{N},\ 0\leq k<2^{i-1}.\] Obviously
$\mathbb{P}(A_{\infty})=1$. Let $p\in(0,1)$ and
$\tau:\mathbb{N}\rightarrow\mathbb{N}$ be any increasing function.
Then
\begin{align*}
\mathbb{E}\Big(\big(\frac{\displaystyle\sum_{i=1}^n\mathbb{I}_{A_{\tau(i)}}}{\displaystyle\sum_{i=1}^n\mathbb{P}({A_{\tau(i)}})}\big)^p\Big)&\leq
\frac{\displaystyle\mathbb{E}
\Big(\big({\sum_{i=1}^n\mathbb{I}_{A_{\tau(i)}}}\big)^p\Big)}{\displaystyle(\frac{n}{2})^p} \ \ \ \ \ \ \  \ \ \ \ \ \ \  \ \ \ \ \ \ \ \ \ \ \ \  \ \ \ \ \ \ \ \ \ (\mathbb{P}(A_{\tau(i)})\geq\frac{1}{2})\\
&\leq\frac{\displaystyle(\int_{[0,\frac{1}{2}]}\sum_{i=1}^n\mathbb{I}_{A_{\tau(i)}}d\mathbb{P})^p\cdot(\frac{1}{2})^{1-p}+\frac{n^p}{2}}
{\displaystyle(\frac{n}{2})^p} \ \ \ \ \ \ (\mbox{H\"{o}lder's inequality})\\
&\leq\frac{\displaystyle(\log_22n)^p\cdot(\frac{1}{2})^{1-p}+\frac{n^p}{2}}
{\displaystyle(\frac{n}{2})^p} \ \
 \big(\mathbb{P}(A_{\tau(i)}\cap[0,\frac{1}{2}])\leq\mathbb{P}(A_{i}\cap[0,\frac{1}{2}])\big).
\end{align*}
Letting first $n\rightarrow\infty$ then $p\rightarrow0$, we get
$\mbox{MS}(\{A_{\tau(n)}\})\leq\frac{1}{2}$. This example shows that
Conjecture 1 is false.

\section{A counterexample to Conjecture 2}

 In the beginning let us  recall two lower
bounds for  $\mathbb{P}(\cup_{i=1}^m A_i)$, where $\{A_i\}_{i=1}^m$
are finitely many events with non-zero probabilities in a
probability space $(\Omega,\mathbb{P})$. First, the Gallot-Kounias
bound (\cite{Gallot,Kounias}, see also \cite{FengLiShen2,Frolov} for
more details) claims that
\begin{equation}\label{GK}
\mathbb{P}(\cup_{i=1}^mA_i)\geq\max_{(\omega_1,\ldots,\omega_m)\in\mathbb{R}^m}
 \frac{\displaystyle\big(\sum_{i=1}^m\omega_i \mathbb{P}(A_i)\big)^2}{\displaystyle\sum_{i=1}^m\sum_{j=1}^m\omega_i\omega_j \mathbb{P}(A_{i}A_j)}
 =\sum_{i=1}^m\gamma_i\doteq \mbox{GK}(\{A_i\}_{i=1}^m),
\end{equation}
 where $\frac{0}{0}\doteq0$ and $(\gamma_1,\ldots,\gamma_m)\in\mathbb{R}^m$ is any solution to
\begin{equation}\label{GK2}
 \Big(\frac{\mathbb{P}(A_{i}A_j)}{\mathbb{P}(A_{i})\mathbb{P}(A_j)}\Big)_{m\times m} \left(
\begin{array}{c}
\gamma_1 \\
\vdots\\
\gamma_m
\end{array} \right )=\left ( \begin{array}{c}
1 \\
\vdots\\
1
\end{array} \right).
\end{equation} Second, Kuai, Alajaji and Takahara (\cite{KAT}, see
also \cite{Hoppe1,Hoppe2,Prekopa}) proved that
\begin{equation}
\label{KAT} \mathbb{P}(\cup_{i=1}^mA_i) \geq \sum _{i=1}^m \Big(
\frac {\displaystyle\theta_i \mathbb{P}(A_i)^2}{\displaystyle S_i+
(1-\theta_i)\mathbb{P}(A_i)}+ \frac
{\displaystyle(1-\theta_i)\mathbb{P}(A_i)^2}{\displaystyle
S_i-\theta_i\mathbb{P}(A_i)} \Big)\doteq
\mbox{KAT}(\{A_i\}_{i=1}^m), \end{equation} where
$S_i\doteq\sum_{j=1}^m\mathbb{P}(A_iA_j)$, $\theta_i$ is the
fractional part of $\frac{S_i}{\mathbb{P}(A_i)}$.

Next let us explain how will we find a  counterexample to Conjecture
2. Suppose  a sequence of events $\{A_n\}_{n=1}^{\infty}$ with
non-zero probabilities occur periodically as follows:
\[A_1,A_2,\ldots,A_m,A_1,A_2,\ldots,A_m,A_1,A_2,\ldots,A_m,\ldots.\]
Obviously,
\[\mathbb{P}(A_{\infty})
\geq\limsup_{n\rightarrow\infty}\mathbb{P}(\bigcup_{i=n}^{n+m-1}A_i)\geq\limsup_{n\rightarrow\infty}\mbox{KAT}(\{A_i\}_{i=n}^{n+m-1})=\mbox{KAT}(\{A_i\}_{i=1}^m).\]
On the other hand, it is easy to observe that
\[\sup\limits_{\tau:\mathbb{N}\rightarrow\mathbb{N}\ \mbox{is increasing}}\mbox{ER}(\{A_{\tau(n)}\})\leq\max_{(\omega_1,\ldots,\omega_m)\in\mathbb{R}^m}
 \frac{\displaystyle\big(\sum_{i=1}^m\omega_i \mathbb{P}(A_i)\big)^2}{\displaystyle\sum_{i=1}^m\sum_{j=1}^m\omega_i\omega_j \mathbb{P}(A_{i}A_j)}
 =\mbox{GK}(\{A_i\}_{i=1}^m).
 \]
Hence to disprove Conjecture 2 it suffices to construct finitely
 many
 $\{A_i\}_{i=1}^m$ so that
 $\mbox{GK}(\{A_i\}_{i=1}^m)<\mbox{KAT}(\{A_i\}_{i=1}^m)$.

  To this aim consider six events $\{A_i\}_{i=1}^6$ in a finite probability space
\begin{eqnarray*}
&\begin{tabular}{|c|c|c|c|c|c|c|c|c|}
\hline $x$ & $\mathbb{P}(\{x\})$ & $A_1$ & $A_2$ & $A_3$ & $A_4$ & $A_5$ & $A_6$\\
\hline $x_1$ & 0.2 & $\bigstar$ & $\bigstar$ & $\bigstar$ & & & $\bigstar$   \\
\hline $x_2$       & 0.2 & $\bigstar$ & & &$\bigstar$ &$\bigstar$ &$\bigstar$ \\
\hline $x_3$       & 0.2 & &$\bigstar$ & $\bigstar$ & & & $\bigstar$ \\
\hline $x_4$       & 0.2 &$\bigstar$ & & & $\bigstar$ & $\bigstar$ &\\
\hline $x_5$       & 0.2 & &$\bigstar$ & $\bigstar$ &$\bigstar$ &$\bigstar$ &  \\
\hline
\end{tabular}
\end{eqnarray*}
with joint probability matrix
\[
(\mathbb{P}(A_iA_j))=\left ( \begin{array}{cccccc}
  0.6  & 0.2 & 0.2 & 0.4 & 0.4& 0.4\\
  0.2 & 0.6 & 0.6 & 0.2 & 0.2& 0.4\\
 0.2 & 0.6 &0.6 & 0.2 & 0.2& 0.4 \\
0.4& 0.2 &0.2 &0.6 & 0.6& 0.2 \\
0.4& 0.2 &0.2 &0.6 &0.6 &0.2\\
0.4& 0.4 &0.4 &0.2 & 0.2&0.6
\end{array} \right ).
\]
Then it is straightforward to work out from
(\ref{GK})$\sim$(\ref{KAT}) that
\[\mbox{GK}(\{A_i\}_{i=1}^6)=\frac{54}{55}<1=\mbox{KAT}(\{A_i\}_{i=1}^6).\]
We are done.


\end{document}